\documentclass[reqno, 12pt]{amsart}

    \usepackage{latexsym,amsfonts, amssymb,amsmath,amscd,amscd,amsthm, mathrsfs,
        amsxtra,epsfig,setspace}

\usepackage{color}

\newtheorem{thm}{Theorem}[section]
\newtheorem{defn}[thm]{Definition}

\newtheorem{rmk}{Remark}

\numberwithin{equation}{section} \numberwithin{thm}{section}
 \numberwithin{rmk}{section}
\numberwithin{figure}{section}

\def\mF{\mathscr{F}}
\def\mG{\mathscr{G}}

    \author{Chunjing Xie}
     \address[Chunjing Xie]{Department of Mathematics, University of Michigan,
              Ann Arbor, MI 48109} 
    \email{ cjxie@umich.edu}
    
    \author{Xiaoming Zheng} 
    \address[Xiaoming Zheng]{Department of Mathematics, Central Michigan University,
              Mount Pleasant, MI 48859 }
 \email{zheng1x@cmich.edu (Corresponding author)}

    \date{}

    \title[Viscoelastic Angiogenesis]{Analysis and simulations of a Viscoelastic Model of Angiogenesis}

\begin{document}
  
    \begin{abstract}
{The work analyzes a one-dimensional viscoelastic model of blood vessel growth under  nonlinear friction with surroundings, and provides numerical simulations for various growing cases. For the nonlinear differential equations, two sufficient conditions are proven to guarantee the global existence of biologically meaningful solutions. Examples with breakdown solutions are captured by numerical approximations. Numerical simulations demonstrate this model can reproduce angiogenesis experiments under various biological conditions including blood vessel extension without proliferation and blood vessel regression.}
    \end{abstract}


    \maketitle
    
\section{Introduction\label{intro}}
Angiogenesis, the growth of new blood vessel from pre-existing vasculature, is crucial to many physiological and pathological processes, such as embryonic development, tumor growth, and wound healing.
Typically these new vessels are very thin (thus also called capillaries) and are lined up by
an array of tightly adhered endothelial cells. Endothelial cells migrate up the gradient of chemotactic cues such as  vascular endothelial growth factor (VEGF) released by tumor cells.
Along a capillary, the tip cell generates a protrusion force to drag the whole capillary forward,
and  at the same time, stalk cells after the tip proliferate and generate new masses to sustain the capillary extension. During the capillary extension, endothelial cells have to overcome the adhesion (or friction) with surroundings, such as other types of cells and extracellular matrix. A detailed description of angiogenesis biology can be found in \cite{Othmer04review}.

Traditional mathematical models on the capillary extension focus on the chemotactic migration of endothelial cells which is typically modeled as a reaction-convection-diffusion equation (see reviews in \cite{Othmer04review} and \cite{JacksonZheng2010}). Such models lack the investigation of  relationships between the growth of endothelial cells,  mechanical extension of the capillary, and  interactions between the capillary and surroundings. These topics are central not only to angiogenesis but also to the biological growth of all soft tissues. However, ``the mathematical modeling of biological growth is currently in an adolescent stage" (\cite{garikipati}), because many constitutive relations are still  unknown. Therefore, it is a big challenge to model the relationships mentioned above.

Among efforts to meet this challenge are one-dimensional models
of  epithelial cell  monolayers migration (\cite{Swigon2007},\cite{Fozard}) in wound healing,
and a one-dimensional viscoelastic angiogenesis model proposed in \cite{ZhengJacksonKoh}. Common features of these models include: 
(i) the blood vessel capillary or epithelial cell mono-layer is regarded as viscoelastic material connected by tight cell-cell junctions;
(ii) the cell growth is treated as an isotropic pressure; 
(iii) the nonlinear friction exists between cells and  surroundings. In this paper, we will first analyze the angiogenesis model in \cite{ZhengJacksonKoh},
and then utilize a numerical scheme to find approximate solutions to this nonlinear problem.

Assume, initially, a blood vessel capillary occupies the entire interval $[0,1]$, and every point of the capillary is denoted by the Lagrangian coordinate, $x$, where $x=0$ denotes the root which is fixed in space, and $x=1$ represents the tip which is free to move. The deformation of the  point $x$ at time $t$ is denoted as $u(x,t)$. The  non-dimensionalized  model is
\begin{equation}
 \beta  \vline 1+\frac {\partial u}{\partial x}  \vline  \frac{\partial u}{\partial t}=  \frac{\partial \sigma}{\partial x}, \label{eq1} 
 \end{equation}
where $\beta$ is the friction between endothelial cells and surroundings, $\sigma$ is the stress tensor defined by
\begin{equation}
 \sigma =  \frac{\partial u}{\partial x} + \mu \frac{\partial \dot{u}} {\partial x }
  - { (f(x,t)-1)}, \label{eq2}
 \end{equation}
$\mu$ is the viscosity of cells, 
$\dot {u}=\frac{\partial u}{\partial t}$ is the material derivative, 
and $f(x,t)$ is the endothelial cell mass density. 
The initial and boundary conditions are
 \begin{eqnarray}
 &&u(x,t=0) = 0 \label{BC2}, \\
&& u(x=0,t) = 0, \quad \sigma (x=1,t) = g(t). \label{BC4}
\end{eqnarray}
where the function $g(t)$ is the protrusion force exerted by the capillary tip cell.

Note that $(1+u_x)$ is the deformation rate of  endothelial cells in the capillary. The equation \eqref{eq1} has the absolute sign of $(1+u_x)$, so it may admit solutions of negative deformation rate. However, a negative deformation rate in one-dimensional space is not biological:  it means one cell simply passes through another cell. Therefore, we only seek biological solutions as defined below.
\begin{defn}
{\bf Biological solution.} If $(1+u_x)\ge 0$, then $u$ is a biological solution.
\end{defn}

To penalize the non-biological solutions, we remove the absolute sign of $(1+u_x)$ in Eq.\,\eqref{eq1}.
For simplicity, we assume  $\beta=1$ and $\mu=0$. The justification to set $\mu=0$ is provided in  \S\,\ref{justification} and the Appendix.
Therefore, the problem we will consider in the main body of this paper becomes
\begin{equation}\label{IBVP}
\left\{
\begin{aligned}
&\left(1+\frac {\partial u}{\partial x}\right)  \frac{\partial u}{\partial t}=  \frac{\partial}{\partial x}\left( \frac{\partial u}{\partial x}-  (f(x,t)-1)\right), \\
& u(0,t) = 0, \,\,\,\, \frac{\partial u}{\partial x}(1, t)= f(1, t)+g(t)-1, \\
&u(x,0) = 0.
\end{aligned}
\right.
\end{equation}
\begin{rmk}
If $(1+u_x)$ happens to be negative, then the first equation of \eqref{IBVP} would behave like a backward heat equation, which is notorious for the ill-posedness of the problem.
\end{rmk}

\begin{rmk}
We will see that the endothelial cell density $f(x,t)$ is dominant in the existence
of biological solutions. Because its presence is  similar to the pressure term in the stress tensor of fluid mechanics, we will also call it the \underline{pressure}.
\end{rmk}

The theoretical result of this paper is stated in Theorem 2.1 of \S\,\ref{secexistence}, which shows that under certain assumptions on $f$ and $g$, the biological solutions exist globally. The key technique to prove Theorem 2.1 is to find some nontrivial transformations to simplify the problem for which the maximum principle can be applied. Because of the nonlinearity, it is hard to find the explicit solutions. Therefore,  in \S\,\ref{secnumerics} we develop a numerical scheme to solve the system \eqref{IBVP} and use it to examine the conditions in Theorem 2.1, including the cases when global biological solutions exist and solutions break down in finite time. In \S\,\ref{biogrowth}, three examples of angiogenesis under different biological conditions are presented using numerical simulations. Finally, we  summarize our analysis in \S\,\ref{conclusion}. An appendix gives some descriptions of the model and the discussion of its parameters.


\section{Analysis of solution existence and asymptotic behavior}\label{secexistence}

Before solving the problem, we investigate the compatibility conditions of the initial  and boundary conditions.  It follows from the initial condition in (\ref{IBVP}) that $u_{xx}=0$ at $(x,t)=(0,0)$. On the other hand, the boundary condition at $x=0$ in (\ref{IBVP}) yields that $u_{t}(0, 0)=0$. Using the differential equation in (\ref{IBVP}), we see that $f$ must satisfy $f_x(0,0)=0$.

\begin{thm}\label{thmexistence}
Suppose $f\in C^2,\, \,g\in C^1$, $f_x(0,0)=0$, $f\geq 0$, $g\geq 0$, and
\begin{equation}\label{condfg}
\inf_{t\geq 0} (f(1,t)+g(t))>0.
\end{equation}
\begin{enumerate}
\item[(a)] Suppose $f_x(0,t)=0$ and $f_{xx}\leq 0$, then there exists a global biological solution for the problem (\ref{IBVP}).
If, moreover, $f(x,t)=f(x)$ and $g(t)=\bar{g}$ (nonnegative constant) when $t\geq T_0$ for some $T_0>0$, then
\begin{equation}\label{asymptotic}
|u(x,t)-\mF(x)-\mG(x)|\leq Ce^{-\alpha t}
\end{equation}
for some $\alpha>0$, where
\begin{equation}
\mF(x)=\int_0^x (f(s)-1)ds,\,\,\mG(x)=\bar{g}x.
\end{equation}
\item[(b)] Suppose $f(x,t)=f(x)$, then there exists a global biological solution for (\ref{IBVP}) provided that
\begin{equation}\label{condoscf}
\max_{x\in [0, 1]} f(x)-\min_{x\in [0, 1]} f(x) < 1.
\end{equation}
If, in addition,  $g(t)=\bar{g}$(constant), then the solution of the problem (\ref{IBVP}) also satisfies (\ref{asymptotic}).
\end{enumerate}
\end{thm}

\begin{rmk}
In part  (a) of the Theorem \ref{thmexistence}, the assumptions  $f_x(0)=0$ and $f_{xx}\leq 0$ indicate that $f$ is a decreasing function. One such an example occurs when the endothelial cells are dying from the tip, so the capillary regresses. In part (b), the assumptions that $f$ is steady in time and its variation is less than unity imply the cell density along the capillary is fixed in time and is quite uniform in space.
\end{rmk}

\noindent{\bf Proof of Theorem \ref{thmexistence}:} Set $F(x,t)=f(x,t)-1$.
\vskip 0.2cm
{\bf (Part a)}. Set
\begin{equation*}
\omega(x,t)=x+u(x,t).
\end{equation*}
Then it follows from \eqref{IBVP} that $\omega$ satisfies the equation
\begin{equation*}
\omega_x\omega_t=\omega_{xx}-F_x(x,t).
\end{equation*}
and
\begin{equation*}
\begin{aligned}
&\omega(0, t)=0,\quad \frac{\partial \omega}{\partial x}(1, t)=f(1, t) +g(t),\\
&\omega(x, 0)=x.
\end{aligned}
\end{equation*}
Therefore,
\begin{equation*}
\omega_t=\frac{\omega_{xx}}{\omega_x}-\frac{F_x(x,t)}{\omega_x}.
\end{equation*}
Because $F_x(0, t)=0$ and  $\omega(0, t)=0$, we have $\omega_{xx}(0, t)=0$. If we extend the functions as
follows
\begin{equation*}
\tilde{\omega}=\left\{
\begin{aligned}
& \omega(x, t),\,\, \text{if}\,\, 0\leq x\leq 1,\\
& -\omega (-x, t),\,\, \text{if}\,\, -1\leq x\leq 0,
\end{aligned}
\right. \,\,  \tilde{F}(x, t)=\left\{
\begin{aligned}
& F(x, t),\,\, \text{if}\,\, 0\leq x\leq 1,\\
& F(-x, t),\,\, \text{if}\,\, -1\leq x\leq 0,
\end{aligned}
\right.
\end{equation*}
then $\tilde{\omega}\in C^2([0, 1])$ and satisfies
\begin{equation}\label{eqomega}
\tilde{\omega}_t=\frac{\tilde{\omega}_{xx}}{\tilde{\omega}_x}
-\frac{\tilde{F}_x(x,t)}{\tilde{\omega}_x}.
\end{equation}
Differentiating equation (\ref{eqomega}) with respect to $x$, one has
\begin{equation*}
(\tilde{\omega}_x)_t=\frac{(\tilde{\omega}_x)_{xx}}{\tilde{\omega}_x}
-\frac{(\tilde{\omega}_{xx})^2}{(\tilde{\omega}_{x})^2}
-\frac{\tilde{F}_{xx}(x,t)}{\tilde{\omega}_x} +
\frac{\tilde{F}_{x}(x,t)\tilde{\omega}_{xx}}{(\tilde{\omega}_x)^2}.
\end{equation*}
Let $\varphi(t,x)=\tilde{\omega}_x$. Then $\varphi$ satisfies 
\begin{equation}\label{eqvarphi}
\varphi_t=\frac{\varphi_{xx}}{\varphi}-\frac{\varphi_x^2}{\varphi^2}+\frac{\tilde{F}_x
\varphi_x}{\varphi^2}-\frac{\tilde{F}_{xx}}{\varphi}.
\end{equation}
If, in addition, $F_{xx}\leq 0$, then
\begin{equation*}
\varphi_t\geq \frac{\varphi_{xx}}{\varphi}-\frac{(\varphi_x-\tilde{F}_x)}{\varphi^2}\varphi_x.
\end{equation*}
By the maximum principle, Lemma 2.3 in \cite{Lieberman}, one has
\begin{equation}\label{maxest}
\min_{\partial' \Omega_T}\varphi\leq  \varphi,\quad \text{for}\,\, (t, x)\in \Omega_T,
\end{equation}
where $\Omega_T= [-1,1]\times [0,T]$ and $\partial'\Omega_T=\partial\Omega_T\backslash \{(x, T)|-1<x<1\}$.
Note that  $f$ and $g$ satisfy (\ref{condfg}), we have
\begin{equation*}
\varphi=\tilde{\omega}_x=1+u_x=1+F+g=f+g\geq \epsilon_1,\,\,\text{at}\,\, x=1,
\end{equation*}
for some $\epsilon_1>0$.
The oddness of the function $\tilde{\omega}$ yields
\begin{equation*}
\varphi(-1,t)=\varphi(1,t).
\end{equation*}
Therefore, the estimate (\ref{maxest}) is equivalent to
\begin{equation}\label{lbdphi}
\varphi\geq \min \{1, \min_{t\geq 0}\{f(1, t) +g(t)\}\}\geq \min\{1, \epsilon_1\}=\epsilon_2>0.
\end{equation}
Define $W=\varphi^2/2$ and
\begin{equation}\label{defV}
h(t)=\max\{\max_{x\in [0,1]} F_{xx}(x,t), 0\}\,\, \,\, \text{and}\,\, \,\,V(t)=\sup_{\partial' \Omega_T}\frac{\varphi^2}{2} +\int_0^t h(s)ds.
\end{equation}
By straightforward computation, the equation \eqref{eqvarphi} can be rewritten as
\begin{equation}\label{eqW}
W_t=\frac{\sqrt{2}}{2}(W^{-1/2}W_x)_x -(\varphi_x -\tilde{F}_x)\frac{1}{2W}W_x -\tilde{F}_{xx}.
\end{equation}
It follows from \eqref{defV} and \eqref{eqW}) that
\begin{equation*}
\left\{
\begin{aligned}
&(W -V)_t\leq  \frac{\sqrt{2}}{2}\Big(W^{-1/2}(W-V)_x\Big)_x -(\varphi_x -\tilde{F}_x)\frac{1}{2W}(W-V)_x,\,\,\text{in}\,\, \Omega,\\
&W-V\leq 0\,\,\,\, \text{on}\,\,\,\, \partial'\Omega,
\end{aligned}
\right.
\end{equation*}
where we used $V_x=0$. 
By the maximum principle again, we have
\begin{equation}\label{upbdphi}
\sup_{\Omega_T}\varphi = \sup_{\Omega_T}\sqrt{2 W} \leq  \sqrt{2V(T)}.
\end{equation}
As long as $\varphi$ is positive and bounded, the equation   (\ref{IBVP}) is uniformly parabolic. Therefore, the global existence of the problem (\ref{IBVP}) is a consequence of Theorem 12.14 in \cite{Lieberman}.

Set  $\Psi=\varphi-\tilde{F}-\bar{g}$.  If $f(x, t)=f(x)$ and $g(x)=\bar g$ for $t\geq T_0$ , then  equation (\ref{eqvarphi}) becomes
\begin{equation}\label{eqpsidiv}
\Psi_t=\left(\frac{\Psi_x}{\varphi}\right)_x
\end{equation}
for $t\geq T_0$.
Note that $\Psi(-1, t)=\Psi(1, t)=0$, multiplying the both side of (\ref{eqpsidiv}) with $\Psi$ and integrating by parts, one has
\begin{equation}\label{intest1}
\partial_t \int_{-1}^1|\Psi|^2(x, t)dx+2\int_{-1}^1\frac{|\Psi_x|^2}{\varphi}(x, t)dx=0.
\end{equation}

By the maximum principle, for $t\geq T_0$,  the solution $\Psi$ of equation (\ref{eqpsidiv}) admits the following estimate
\begin{equation*}
\Psi\leq \max \{\sup_{-1\leq x\leq 1}\Psi(x, T_0), 0\}.
\end{equation*}
This implies
\begin{equation*}
\sup_{\substack{-1\leq x\leq 1,\\ t\geq T_0}} \varphi \leq  \sup_{-1\leq x\leq 1}(\tilde{F} +\bar{g})+ \max \{\sup_{-1\leq x\leq 1}\Psi(x, T_0), 0\}.
\end{equation*}
Combining with \eqref{upbdphi} gives 
\begin{equation}\label{ubdphi}
\sup_{\substack{-1\leq x\leq 1,\\ t\geq T_0}} \varphi \leq \frac{1}{2\alpha_1}
\end{equation}
for some constant $\alpha_1>0$.
Thus, taking \eqref{lbdphi} and \eqref{ubdphi} into account, it follows from \eqref{intest1} that we have
\begin{equation*}
\partial_t \int_{-1}^1|\Psi|^2(x, t)dx+\alpha_1 \int_{-1}^1|\Psi_x|^2(x, t)dx\leq 0.
\end{equation*}
Using Poincare inequality gives
\begin{equation*}
\partial_t \int_{-1}^1|\Psi|^2(x, t)dx+ 2\alpha \int_{-1}^1|\Psi|^2(x, t)dx\leq 0
\end{equation*}
for some $\alpha>0$. This implies
\begin{equation}\label{energydecay}
\int_{-1}^1 |\Psi(x, t)|^2 dx\leq Ce^{- 2\alpha t}.
\end{equation}
Furthermore, the  Nash-Moser iteration (Theorem 6.17 and Theorem 6.30 in \cite{Lieberman}) shows that that there exist $t_0$ and $C$ (independent of $t$) such that
\begin{equation*}
\sup_{-1\leq x\leq 1}|\Psi(x, t)|\leq C \left(\int_{t-t_0}^{t+t_0} \int_{-1}^1 |\Psi(x,\tau)|^2 dx d\tau\right)^{1/2}. 
\end{equation*}
Using estimate (\ref{energydecay}) yields that
\begin{equation*}
\sup_{-1\leq x\leq 1}|\Psi(x, t)|\leq C e^{-\alpha t}.
\end{equation*}
This finishes the proof for part (a). \\

\vskip 0.5cm
{\bf (Part b)}. Let us start with the equation (\ref{eqvarphi}). If
$F(x, t)=f(x)$ and set
\begin{equation*}
\psi(x,t)=\varphi(x,t)-\tilde{F},
\end{equation*}
then $\psi$ satisfies the equation
\begin{equation*}
\psi_t=\frac{1}{\varphi}\psi_{xx}-\frac{\varphi_x}{\varphi^2}\psi_x.
\end{equation*}
By the maximum principle, Lemma 2.3 in \cite{Lieberman}, one has
\begin{equation*}
\min_{\partial'\Omega}\psi\leq\psi\leq \max_{\partial'\Omega}\psi
\end{equation*}
Thus
\begin{equation*}
\min_{\partial'\Omega} \psi +\tilde{F}\leq\varphi \leq
\max_{\partial'\Omega} \psi +\tilde{F}.
\end{equation*}
Note that at $x=1$,
\begin{equation*}
\psi=\varphi-\tilde{F}=1+F+g-F=1+g,
\end{equation*}
so
\begin{equation*}
\min\psi(1,t)+\tilde{F}(1)=1+\min g(t)+f(1)-1=f(1)+\min g.
\end{equation*}
When $t=0$,
\begin{equation*}
\min(\psi+\tilde{F})=\min(1-( f-1))+(f-1)=1+(f-\max f).
\end{equation*}
If (\ref{condoscf}) holds, then there exists an $\epsilon>0$ such that
\begin{equation*}
\max f -\min f<1-\epsilon.
\end{equation*}
Hence $\varphi>\epsilon$. Similarly, 
\begin{equation*}
\max_{\partial'\Omega}\psi +\tilde{F}= \max\{f(1)+\max g, 1+f-\min f\}.
\end{equation*}
Hence, we can also show global existence of the problem (\ref{IBVP}) via Theorem 12.14 in \cite{Lieberman}.

Similar to part (a), we can get the asymptotic behavior of the solution by energy estimate and Nash-Moser iteration. This finishes the proof for part (b).
\qed

\section{Numerical studies of existence of biological solutions\label{secnumerics}}
The assumptions on $f$ and $g$ in Theorem \ref{thmexistence} are sufficient conditions to guarantee global existence of solutions, thus it would be interesting to examine the solutions when these assumptions are violated. In this section we utilize a numerical method to find approximate solutions of the nonlinear problem. 

\subsection{Numerical method\label{scheme}} We consider solving a more general nonlinear problem
\begin{equation}\label{IBVP_viscous}
\left\{
\begin{aligned}
&\beta \left(1+\frac {\partial u}{\partial x}\right)  \frac{\partial u}{\partial t}=  \frac{\partial}{\partial x}\left( \frac{\partial u}{\partial x} + \mu \frac{\partial \dot{u}}{\partial x}
-  (f(x,t)-1)\right), \\
& u(0,t) = 0, \,\,\,\, \frac{\partial u}{\partial x} + \mu \frac{\partial \dot{u}} {\partial x }
  - { (f(x,t)-1)} = g(t), \\
&u(x,0) = 0.
\end{aligned}
\right.
\end{equation}
with a finite element method in the Sobolev space $H^1(0,1)$, that is, the square-integrable functions up to the first order weak derivative. 

Choose a uniform time step $\Delta t > 0$ and denote the time points when solutions are sought  as $t^k=k\Delta t, k=0, 1, \cdots$. The numerical scheme is to find $u^{k+1}(x) \in H^1(0,1)$ with $u^{k+1}(0)=0$, such that for any test function $\phi\in H^1(0,1)$ with $\phi(0)=0$,
\begin{equation}
 \int\limits_0^1 \beta  (1+\frac {\partial u^{k}}{\partial x})   \frac  {u^{k+1}-u^k} {\Delta t}     \phi =
 \phi(1) g(t^{k+1}) - \int\limits_0^1
 \left(\frac{\partial u^{k+1}}{\partial x} + \mu \frac{\partial} {\partial x }  \frac  {u^{k+1}-u^k} {\Delta t}  - { (f^{k+1}-1)}\right) \phi_x.
\label{peipei2}
 \end{equation}
The spatial domain $[0,1]$ is uniformly discretized into $n$ equal-sized sub-intervals with mesh size $h=\frac 1 n$, and mesh points are denoted as $x_i=(i-1)h, i=1, \cdots, n+1$. The space $H^1(0,1)$ is approximated by  $P^1$ finite element space:
$$
V_h(0,1) = \{ v_h\in C^0(0,1): v_h \text{ is a linear function on each subinterval } [x_i,x_{i+1}], i=1,\cdots, n   \}.
$$
Numerical tests show this scheme is first order accurate in time (data not shown). A Matlab version of the code has been provided in the website http://www.cst.cmich.edu/users/zheng1x/. In all the numerical simulations in this work, we have chosen $n=1000$ and $\Delta t=10^{-5}$, and each simulation result is almost identical to that with the more refined choice  $n=2000$ and $\Delta t=5\times10^{-6}$.
\subsection{Comparisons between  friction and viscosity\label{justification}}
The analysis in Appendix \ref{secapp} suggests that the viscosity and friction have an additive effect on the time scale of capillary extension, and the viscosity of  biological values is negligible compared with the friction.  In this section, we use numerical simulations to test the effect of viscosity on solutions with respect to friction.

\begin{figure}[htbp]
\begin{center}
\psfig{file=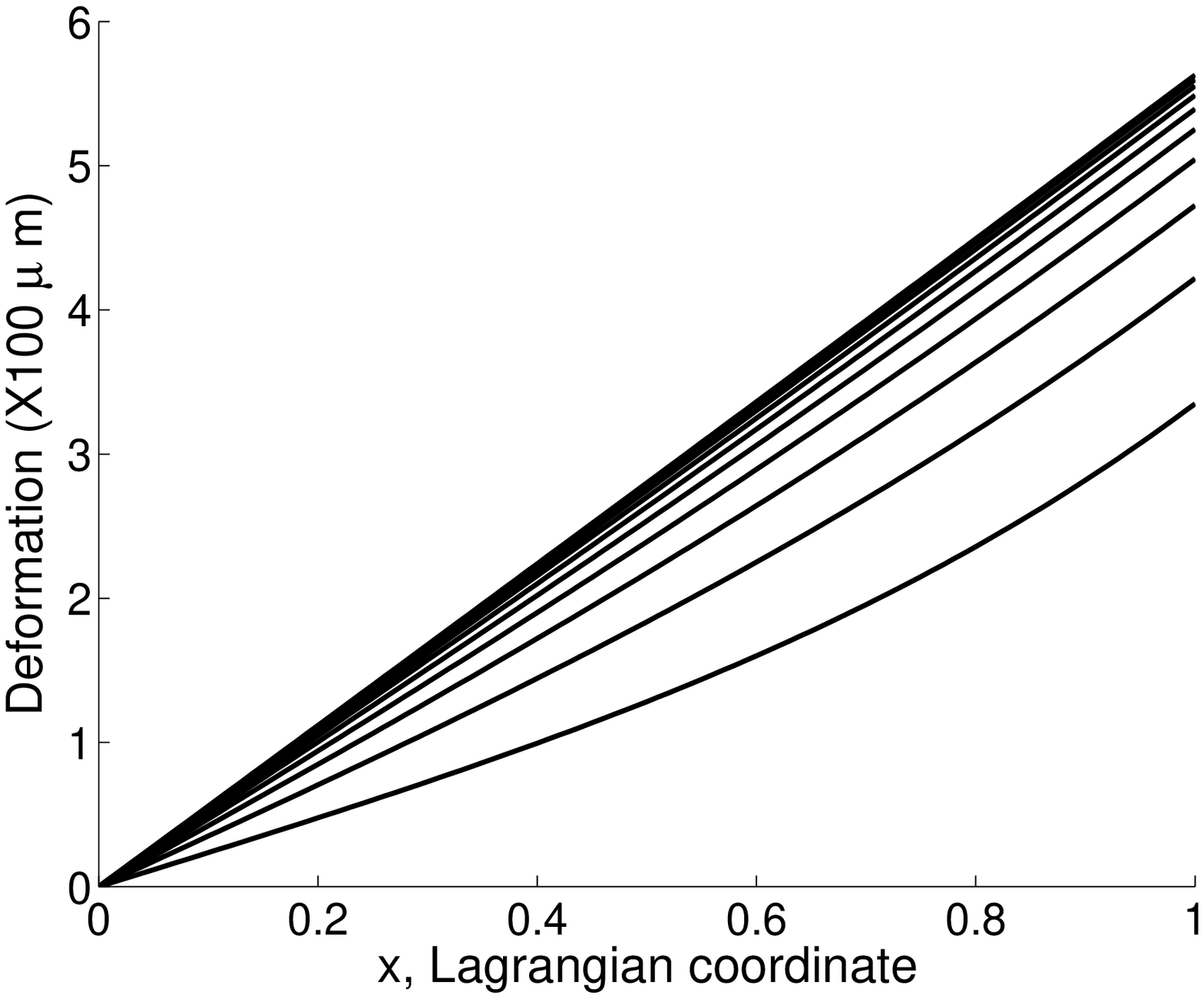, width=2.5in}[a]
\psfig{file=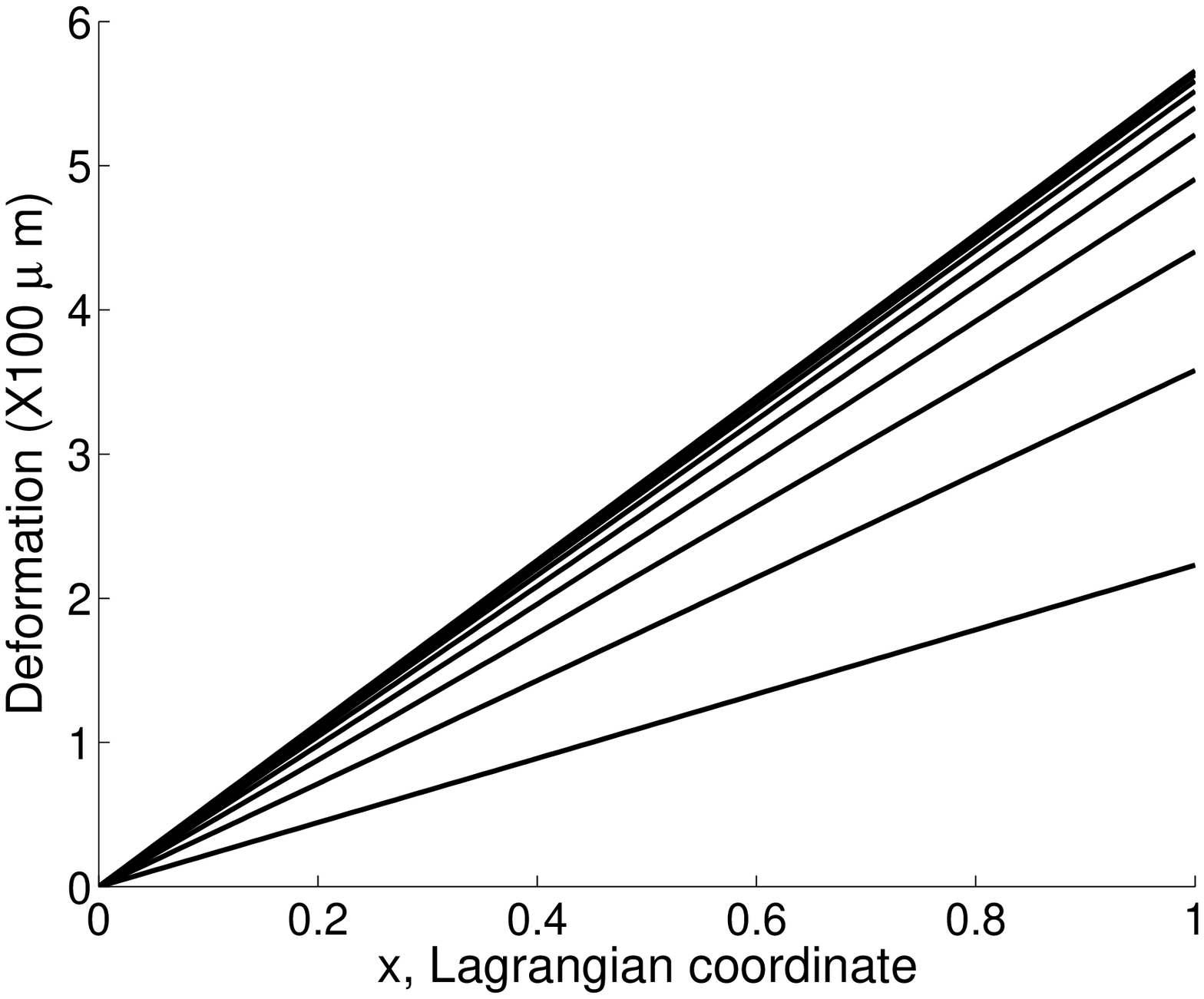, width=2.5in}[b]
\caption{Solutions when $\beta=0.01, g=5.7, f=1$.
[a]: biological case, $\mu=10^{-4}$. The solution curves are at time =$0.01, 0.02, \cdots, 0.1$ (bottom to top). [b]: exaggerated case, $\mu=1$. The solution curves are at time=$0.5, 1, \cdots, 5$ (bottom to top).}
\label{comp1}
\end{center}
\end{figure}

First, we fix $\beta=0.01$, $g=5.7$, $f=1$, and compare the solutions in three cases: inviscid case, $\mu=0$; biological case, $\mu=10^{-4}$ (corresponding to dimensional viscosity $\tilde\mu=10^4\,\frac{pN \cdot s}{\mu m^2}$); and exaggerated case, $\mu=1$ (corresponding to biological viscosity $\tilde\mu=10^8\,\frac{pN \cdot s}{\mu m^2}$, comparable to that of pitch at $20^\circ C$). The numerical results are shown in Fig.\,\ref{comp1}. The $\mu=0$ (results not shown) is almost identical to the $\mu=10^{-4}$ case and  in both cases  the solutions reach the same steady state at roughly the same time $t=0.1$. In the exaggerated case, the solution reaches the steady state at $t=5$. Therefore, when the viscosity is significantly larger than friction, the solution goes to the steady state much slower.

\begin{figure}[htbp]
\begin{center}
\psfig{file=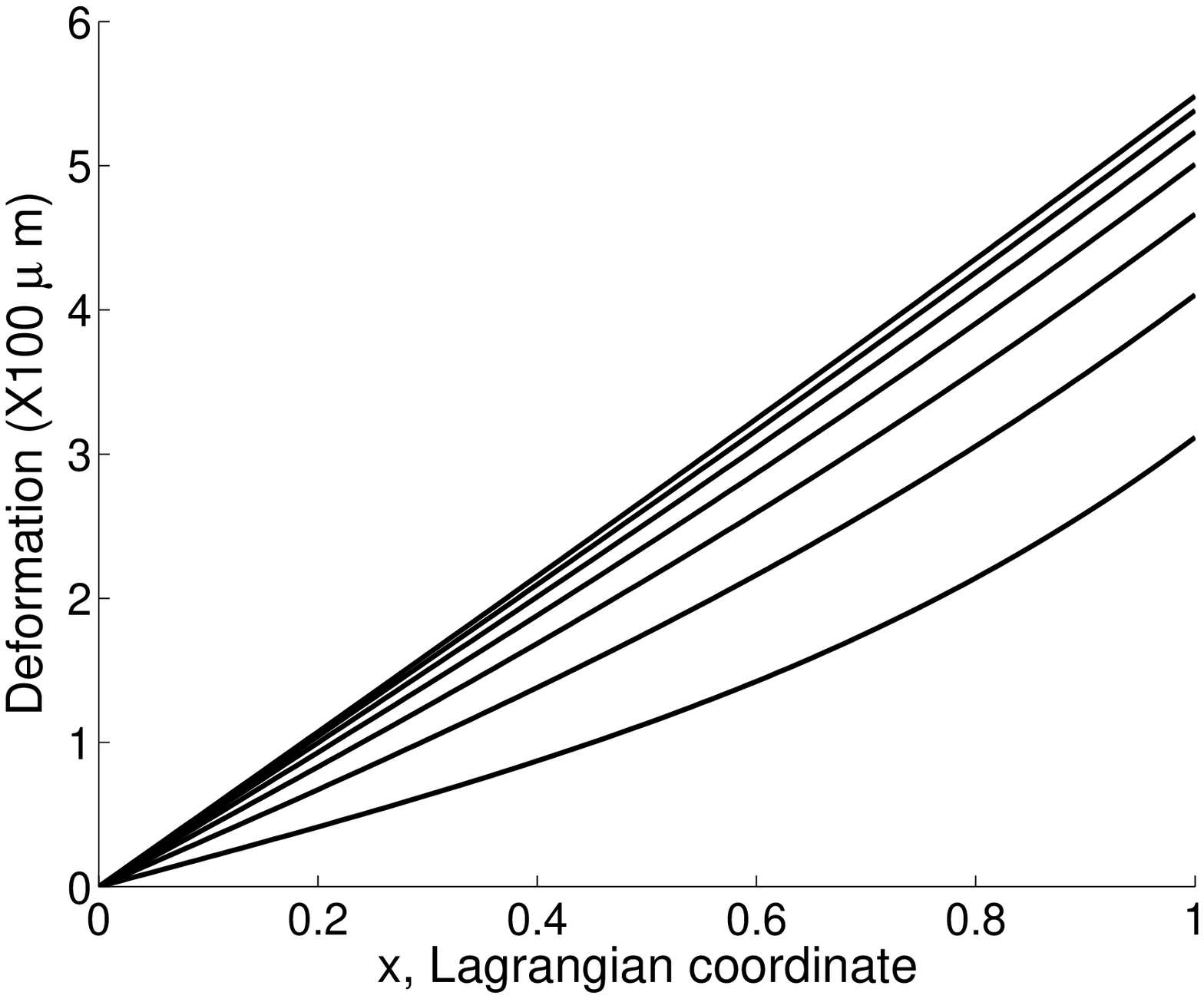, width=2.5in}[a]
\psfig{file=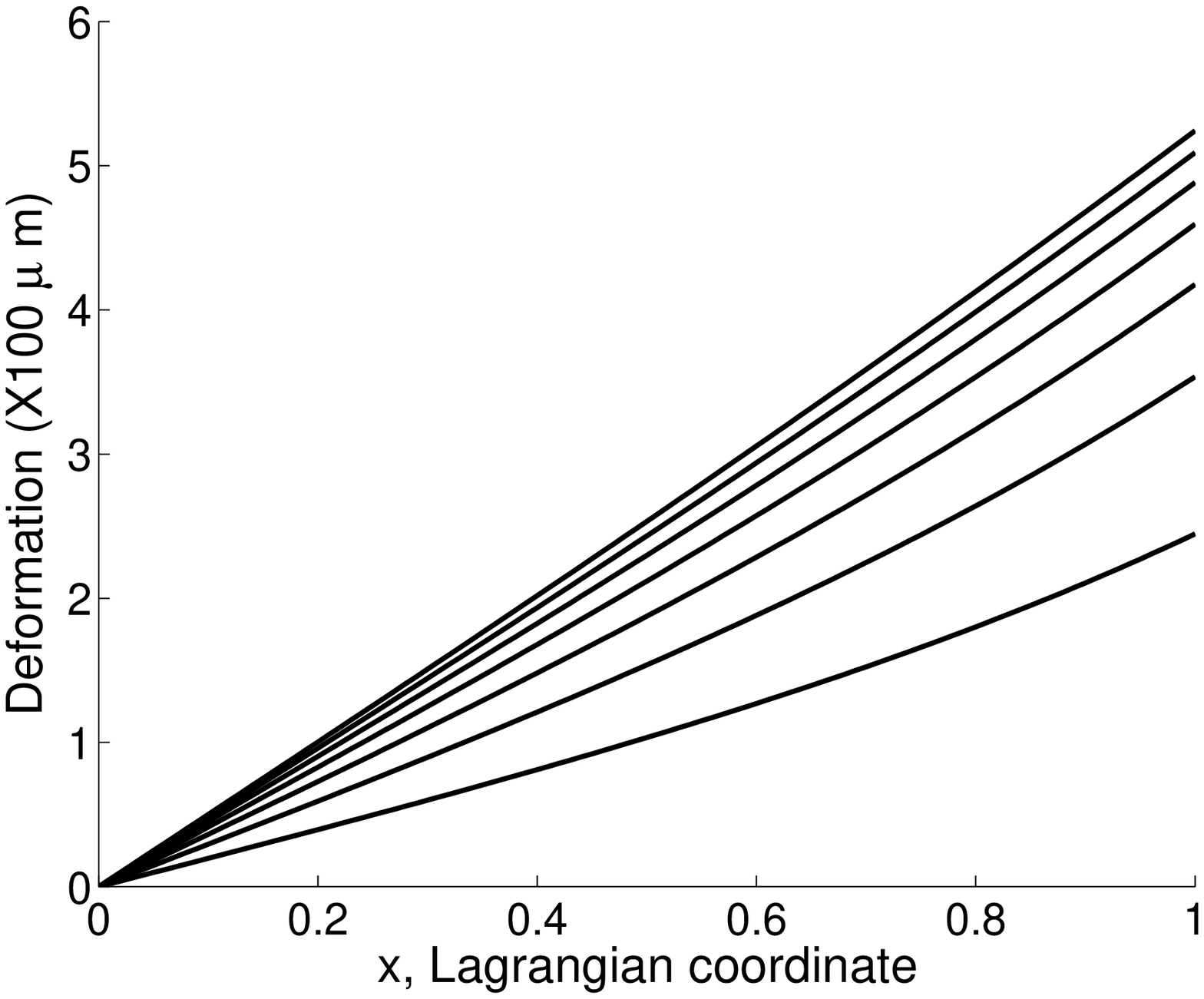, width=2.5in}[b]
\caption{Solutions when $\beta=1$, $g=5.7$, and $f=1$. In both [a] and [b], the solution curves are at $t=1,2,\cdots, 7$ (bottom to top). [a]: biological case, $\mu=10^{-4}$. [b]: exaggerated case, $\mu=1$.}
\label{comp2}
\end{center}
\end{figure}

Second, we change the friction to $\beta=1$ and re-do the above three simulations.  The numerical results are shown in Fig.\,\ref{comp2}. The solutions in the inviscid and biological cases are also not distinguishable, and both reach the steady state at $t=7$, while the solution in the exaggerated case gets to the steady state at $t=15$. Although the $\mu=1$ case has the decreased capillary extension speed,  it is not significantly different from the inviscid case (Fig.\,\ref{comp2}[b]).

These comparisons confirm the suggestion at the beginning of this section, that is, when we focus on the biological values, it is valid to  neglect the viscosity.

\subsection{Examples examining Theorem\,\ref{thmexistence} (a)\label{parta_examples}}
In all the examples of this subsection, we let $\beta=1$, $\mu=0$, and $g=0$.

\begin{figure}[htbp]
\begin{center}
\psfig{file=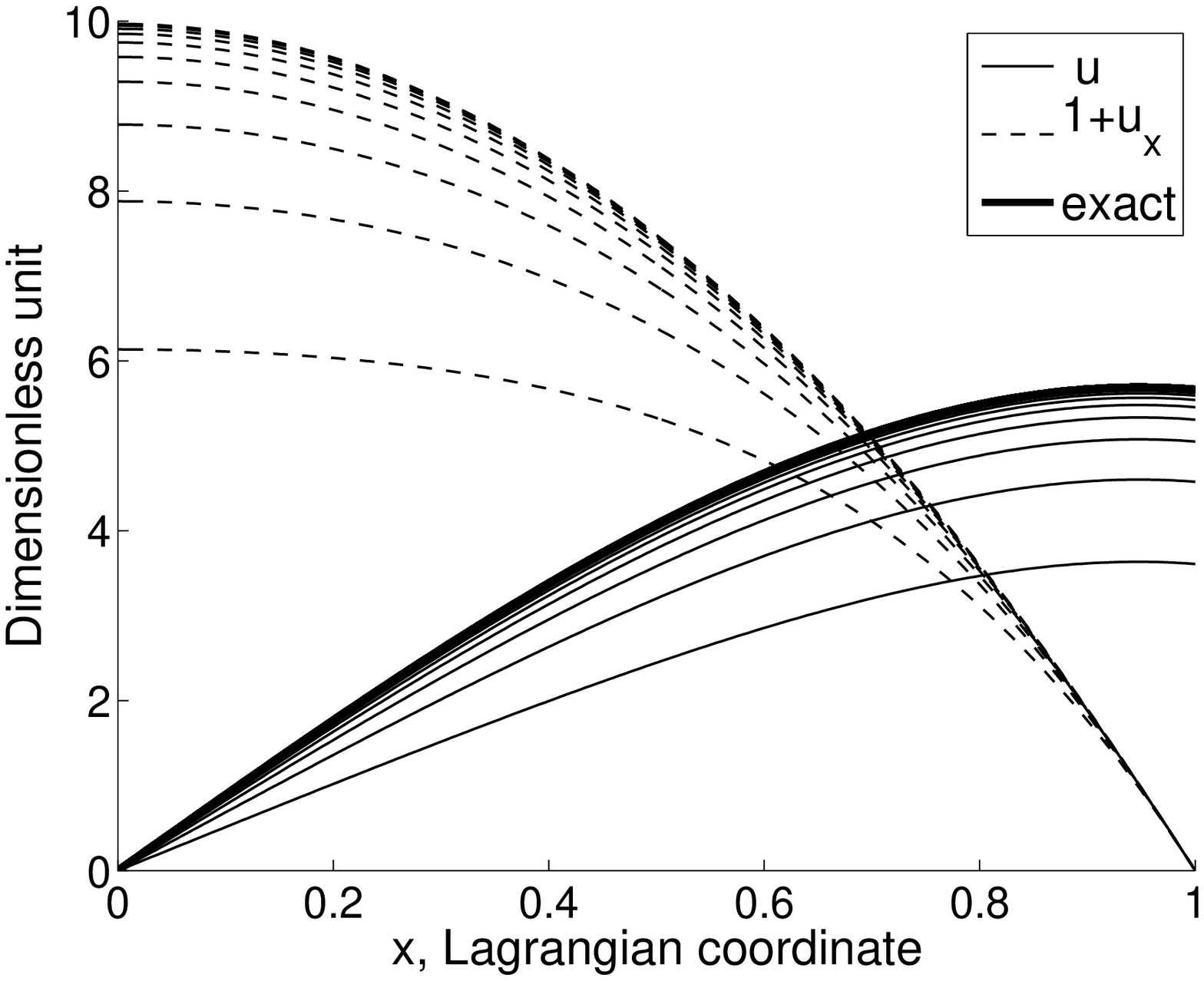, width=2.5in}[a]
\psfig{file=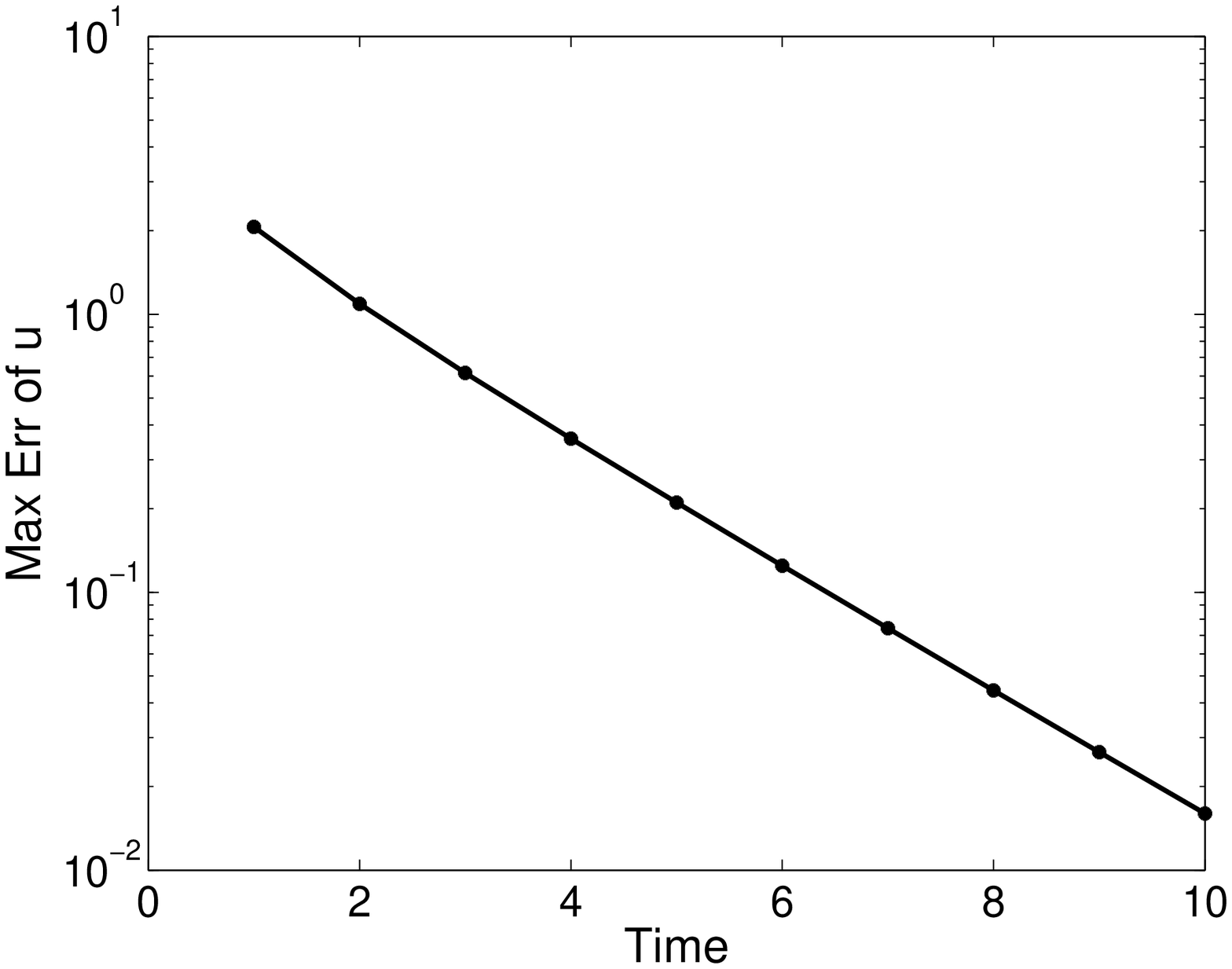, width=2.5in}[b]
\caption{[a] Solution curves (solid) and $(1+u_x)$ curves (dashed) in the case of $\beta=1$, $\mu=0$, $g=0$, $f=-10x^2+10+10^{-6}$, at time $t=1,2,\cdots, 10$ (both from lower to upper). [b] $|u_h-u_s|_\infty$ with respect to time. The y-axis is in log scale.}
\label{example_qq}
\end{center}
\end{figure}

First, we choose  $f=-10x^2+10+10^{-6}$. This choice satisfies the conditions
in Theorem \,\ref{thmexistence} (a). The numerical solution from $t=1$ to $t=10$ is illustrated in Fig.\,\ref{example_qq}[a]. It is clear that  $1+u_x$ is non-negative for all time, and
the solution converges to $u_s=-\frac{10}{3} x^3 + (9+10^{-6})x$, the steady state solution. The maximum error between the numerical solution $u_h$ and $u_s$ is shown in Fig.\,\ref{example_qq}[b], and it tends to zero
in time exponentially.

\begin{figure}[htbp]
\begin{center}
\psfig{file=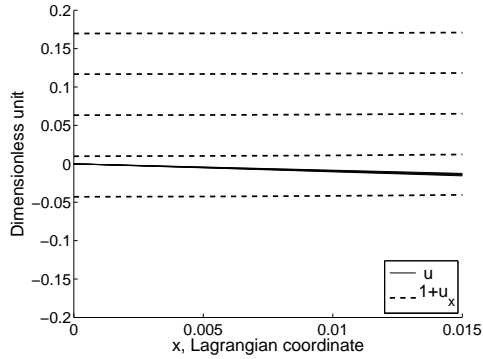, width=2.5in}
\caption{Numerical solution and $1+u_x$ at $t=0.026, 0.027, \cdots, 0.03$ (upper to lower) in the case of $\beta=1$, $\mu=0$, $g=0$, $f=10x^2+10$. }
\label{example_pop}
\end{center}
\end{figure}

Second, we choose  $f=10x^2+10$. The choice of $f$ violates the condition $f_{xx}\le 0$
in Theorem\,\ref{thmexistence} (a). The value of $1+u_x$ is negative for $x\in [0,0.015]$ when $t=0.03$ (Fig.\,\ref{example_pop}). A few time steps later the numerical solution blows up near $x=0$. 

Notice in the mathematical model the density serves as a pressure exerting in the direction $-\nabla f$, that is, from high pressure regions to low pressure regions. The breakdown of the second example stems from the high density at the tip driving towards the root. At the root, the boundary condition $u=0$ makes the root cell  nowhere to escape.  So even after it is smashed to zero length, the pressure still cannot be balanced, which results in the breakdown of the model.  In the first example, the pressure with the same magnitude of gradient exerts on the system but towards the tip direction. The free boundary condition at the tip releases the pressure and avoid the breakdown of the system.

 These two examples implies that, if the pressure is monotonic along the vessel,  the position of high pressure relative to the fixed boundary 
is important for the existence of biological solutions. But positioning the high pressure  at the tip does not necessarily break down the system. For instance, we also tested $f=1+x^2$, in which case the high pressure is at the tip but the global biological solution still exists. Comparing with the case of $f=10+10x^2$, which is of higher pressure drop from the tip to the root, we deduce that the pressure drop  from tip to root must be large enough to break down the system.

\subsection{Examples examining Theorem \,\ref{thmexistence} (b)\label{partb_examples}}
The same as the previous section, we also set $\beta=1$, $\mu=0$, and $g=0$ for all the examples in this subsection.

\begin{figure}[htbp]
\begin{center}
\psfig{file=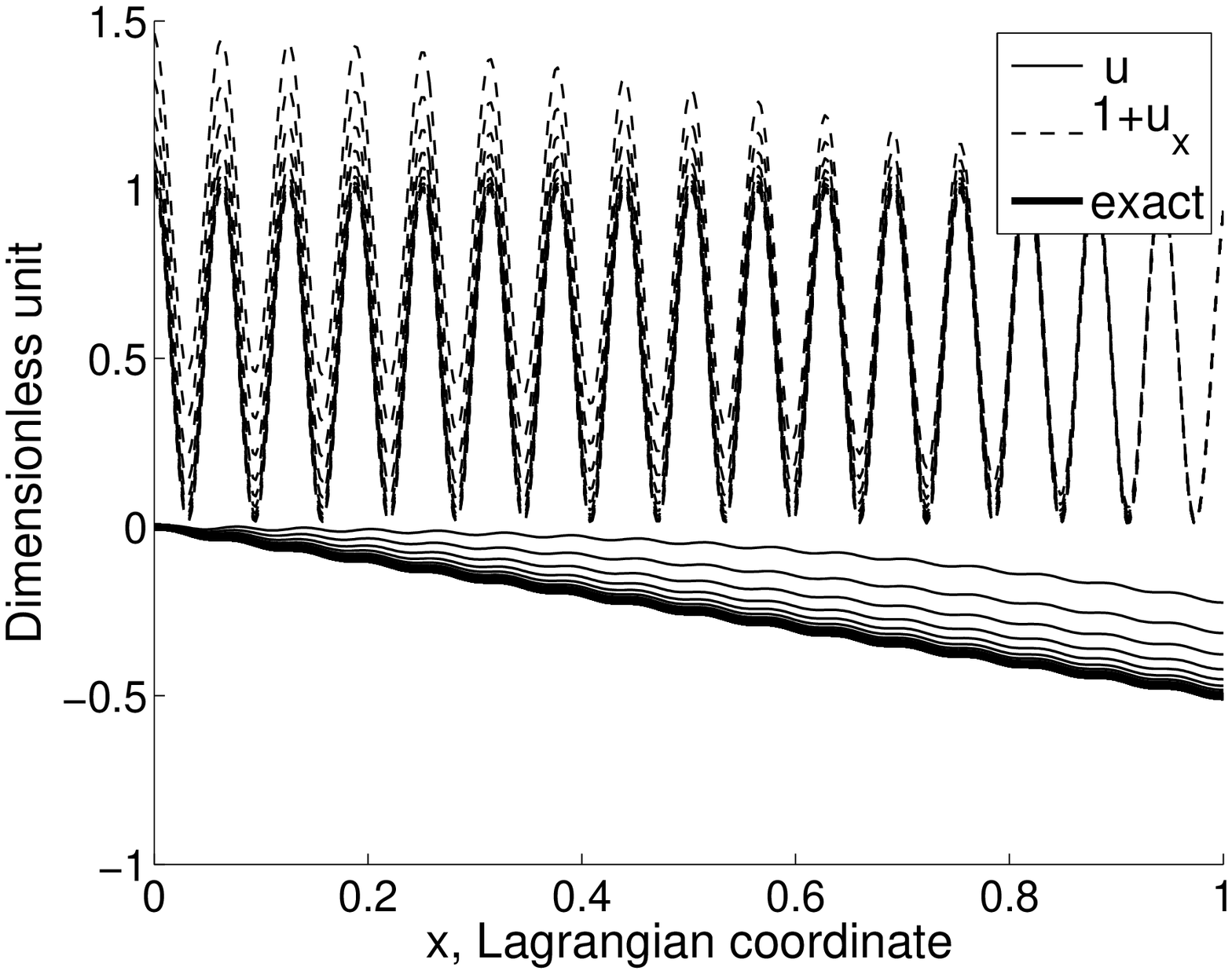, width=2.5in}[a]
\psfig{file=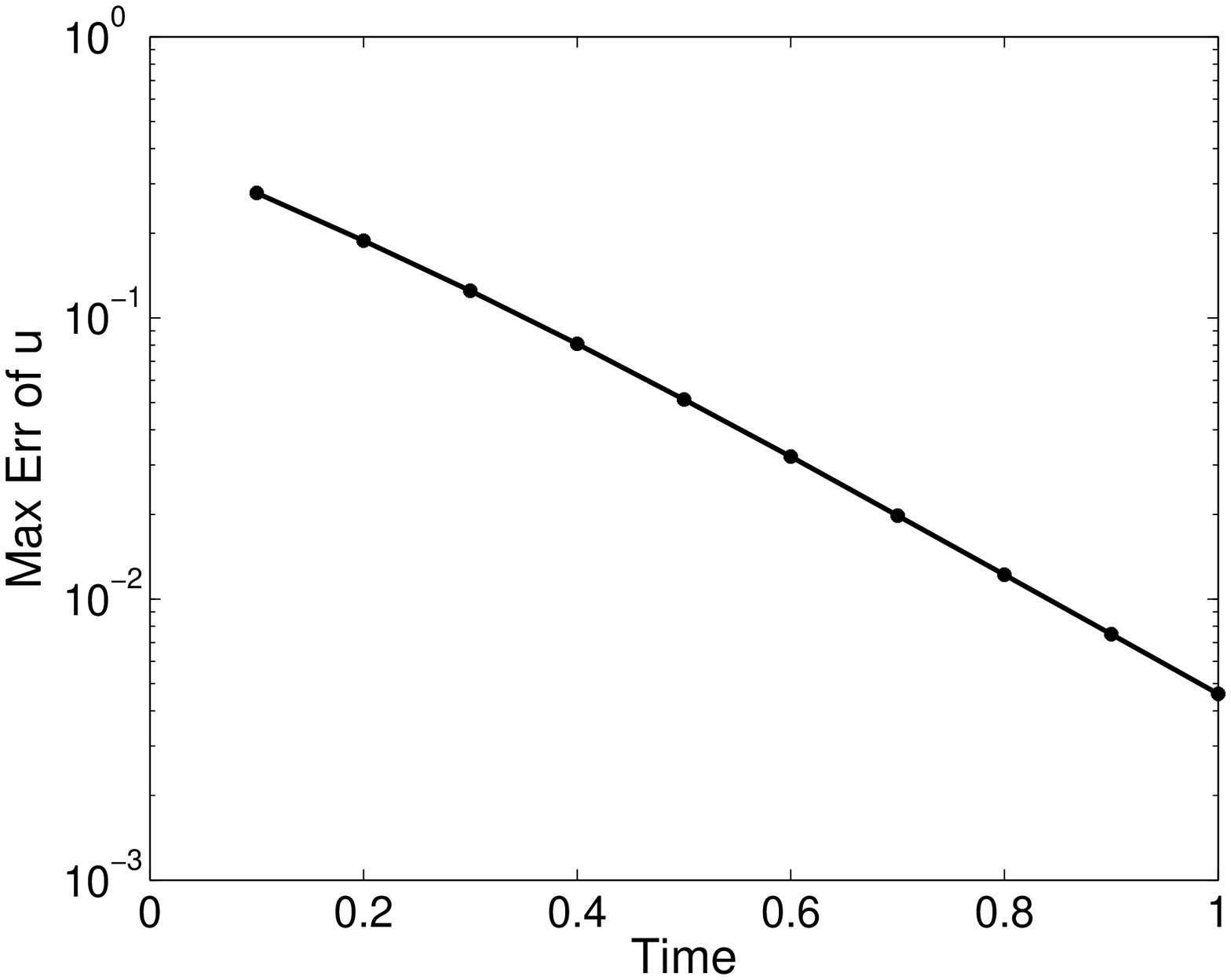, width=2.5in}[b]
\caption{[a] Solution curves (solid) and $(1+u_x)$ curves (dashed) in the case of $\beta=1$, $\mu=0$, $g=0$, $f=0.4999\cos(100x)+0.5$, at time $t=.1,.2,\cdots, 1$ (from upper to lower).
[b] $|u_h-u_s|_\infty$ with respect to time. The y-axis is in log scale.}
\label{example_gg}
\end{center}
\end{figure}

First, let $f=0.4999\cos(100x)+0.5$. This choice satisfies the conditions in Theorem\,\ref{thmexistence} (b). The numerical solution from $t=0.1$ to $t=1$ (Fig.\,\ref{example_gg}[a]) indicates that  $(1+u_x)$ is non-negative for all time, and the solution converges to $0.004999\sin(100x) -0.5x$, the steady state solution. The maximum error between the numerical solution and the steady state solution is shown in Fig.\,\ref{example_gg}[b], and it tends to zero in time exponentially.

\begin{figure}[htbp]
\begin{center}
\psfig{file=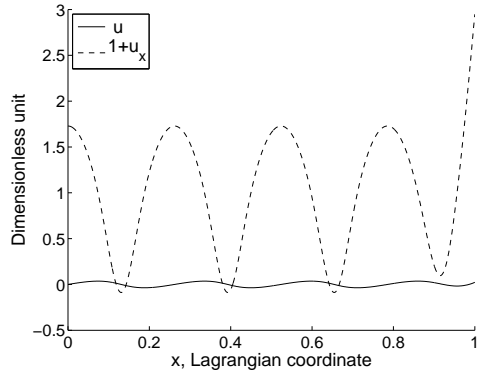, width=2.5in}
\caption{ Solution curve (solid) and $(1+u_x)$ curve (dashed) in the case of $\beta=1$, $\mu=0$, $g=0$, $f=2\cos(24x)+2.1$, at time $t=10^{-3}$.}
\label{example_rr}
\end{center}
\end{figure}

Second, we change to  $f=2\cos(24x)+2.1$. This choice violates the condition $\max_{x\in [0, 1]} f(x)-\min_{x\in [0, 1]} f(x) < 1$ in Theorem \,\ref{thmexistence} (b). The numerical solution at $t=10^{-3}$ is illustrated in Fig.\,\ref{example_rr}, which shows that  $1+u_x$ is negative at several points, where the solution blows up later. 

Note that the pressure gradient in the blow-up case is even smaller than that in the non-blow-up case, therefore, the blow-up in the second example cannot be attributed to the high pressure gradient. However, the pressure amplitude in the blow-up case is larger than the non-blow-up case. Therefore, it is pressure's amplitude that breaks down the system.

\section{Numerical simulations of angiogenesis\label{biogrowth}}
In this section, we apply the mathematical model to simulate  angiogenesis under three real biological conditions: capillary extension without proliferation, extension with proliferation, and capillary regression.

\subsection{Capillary extension without proliferation, but with variable friction\label{ang_1}}
\begin{figure}[htbp]
\begin{center}
\psfig{file=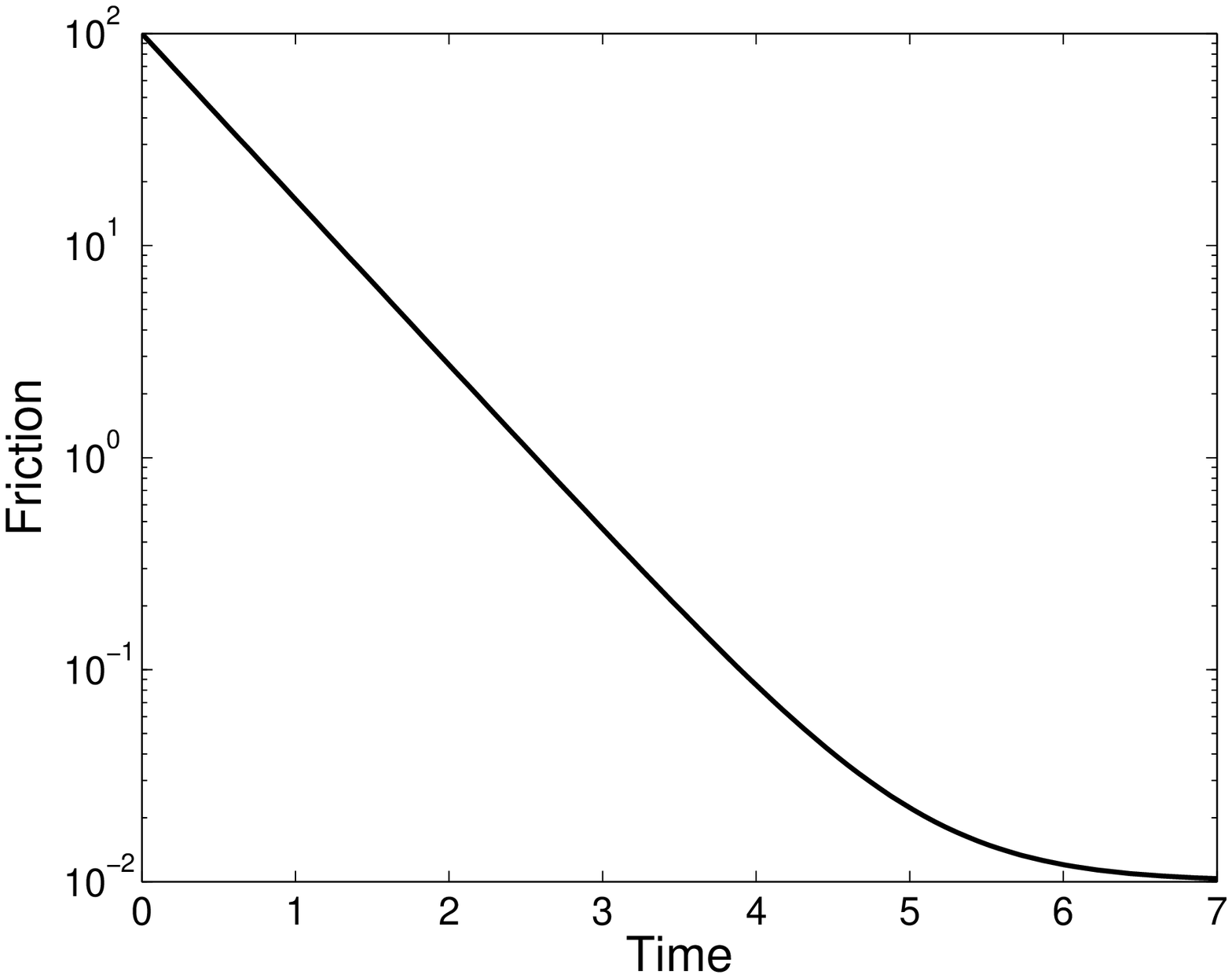, width=2.5in}[a]
\psfig{file=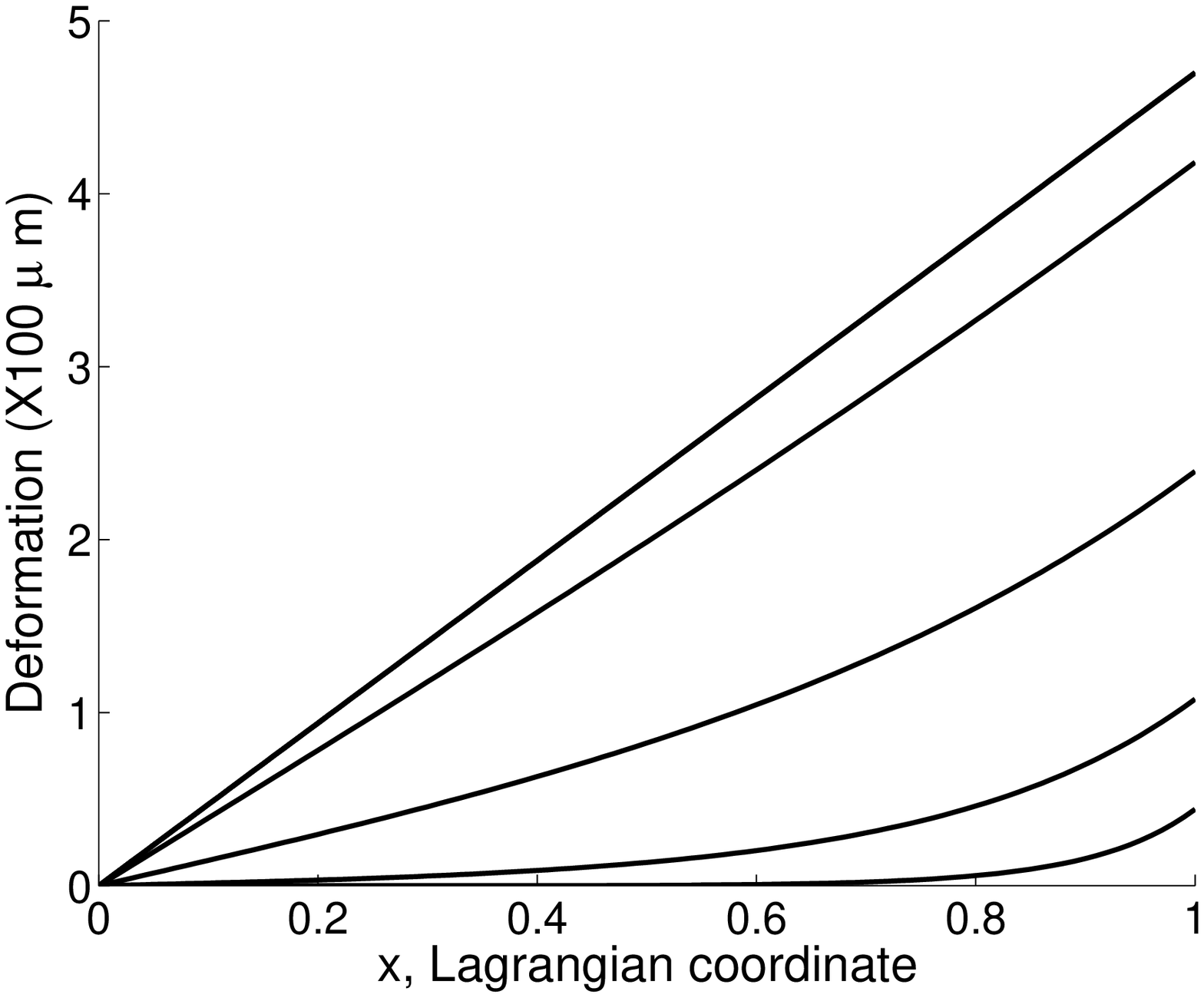, width=2.5in}[b]
\caption{Extension of capillary without proliferation, $f(x,t)=1$. [a] Friction as a function of time, $\beta=0.01+100 e^{-1.8t}$. Y-axis is in log scale. [b] Deformation curves at time $t=1,2,\cdots, 7$ (lower to upper). The curves of $t=5,6,7$ almost coincide. Other parameters: $g=4.7$, $\mu=10^{-4}$. }
\label{fig1}
\end{center}
\end{figure}

In the experiments of \cite{Sholley}, the endothelial cells are given X-ray irradiation. With enough doses of irradiation,  DNA synthesis are stalled and cells lose the ability to proliferate. This means the cell density will keep the same from the beginning. Therefore, the cell density is always equal to the reference density, that is, $f(x,t)=1$ for all the points and all the time (see Appendix for non-dimensionalization). 
The viscosity and protrusion force are chosen as $\mu = 10^{-4}$ and $g = 4.7$, within the biological ranges as discussed in the Appendix. However, the friction is a variable in time. Because, at the initial stage of angiogenesis, endothelial cells have to overcome the tight association with surrounding cells and basement membrane.These obstacles are rapidly reduced after the initial stage, so the friction is smaller with the progression of angiogenesis.  Therefore, we assume the friction is a time-dependent function 
\begin{equation}
\beta(t) =  0.01 + 100e^{-1.8t},
\end{equation}
whose graph is plotted in Fig.\,\ref{fig1}[a]. This is a pretty rough approximation of the true variation of the friction. Indeed the friction depend on many other conditions, such as the proteins regulating EC-environment adhesion, association with pericytes, etc. Readers are referred to \cite{ZhengJacksonKoh} for more details.

The numerical solutions are shown in Fig.\,\ref{fig1}[b]. At time $t=1$, only the points near the tip  have observable motion, because the friction is very large initially and only the cells near the tip can be dragged with remarkable distance by the protrusion force exerted at the tip. Up to $t=3$, all points migrate toward the tip direction. The deformation of the tip is $3.9$ at $t=4$ and $4.7$ at $t=7$. These correspond to the whole capillary length $490\,\mu m$ at day 4 and $570\,\mu m$ at day 7 after the initiation of angiogenesis, very close to the rat corneal experimental results in \cite{Sholley}. 

Although there is no cell proliferation, the protrusion force of the tip cell drags the capillary to extend.  The ultimate deformation is limited by the elasticity and is given by the steady state solution $u_s = g/E$. For the above example, this value is $4.7$, which has been achieved by the numerical simulation. However, the lack of proliferation does affect the capillary extension: the capillary only reaches $570\,\mu m$, but the distance from the root to the chemotactic source is about $2\, mm$ in the experiments of \cite{Sholley}. 
The limited capillary extension without proliferation is an important feature of biological growth, but it is very hard to capture using reaction-diffusion-convection models. For example, 
many reaction-convection-diffusion models (e.g. \cite{Anderson}) only tracks the motion of the tip cell, which is modeled as migrating towards the chemotactic source without any restriction of trailing cells. Thus, it falsely predicts the capillary reaches the chemotactic source even when there is no cell proliferation.
\subsection{Capillary extension with proliferation and variable friction\label{ang_2}}
\begin{figure}[htbp]
\begin{center}
\psfig{file=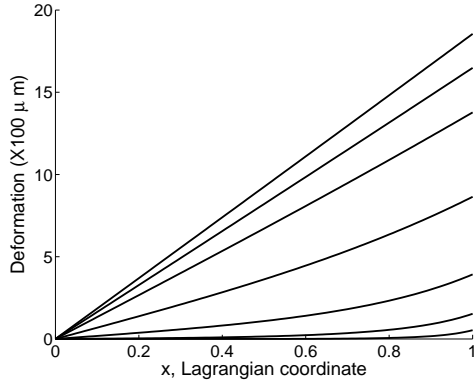, width=2.5in}
\caption{Capillary extension with cell density $f=2t+1$. Other parameters: $\beta=0.01+100 e^{-1.8t}$, $g=4.7$, $\mu=10^{-4}$. Curves from lower to upper: $t=1,2,\cdots, 7$ days.}
\label{fig2}
\end{center}
\end{figure}

Comparing with the last example, we change the cell density to $f=2t+1$. The biological meaning of this function is that the cells uniformly proliferate in the whole capillary: each cell generates two more cells per day. The numerical results are shown in Fig.\,\ref{fig2}.
The value of $u$ of the capillary tip is equal to $0.53$ at $t=1$, $8.64$ at $t=4$, and $18.5$ at $t=7$. These represents the capillary reaches $63\,\mu m$ at Day 1, $964\,\mu m$ at Day 4, and $1950\,\mu m$ at Day 7. These results reproduce the data in the experiments of \cite{Sholley} and at Day 7 the capillary has already reached the chemotactic source.

In this example, both the tip cell protrusion and stalk cell proliferation contribute to the capillary extension. Comparing with the last example, cell proliferation is more decisive in supporting the extension of the capillary until it reaches the chemotactic source. Mathematically, the increase of cell density alters the balance between the previous stress and density. The cells enlarge the deformation to adjust to the higher density, which allows the extension of the vessel. Therefore, the tip cell migration plays a leading role in capillary extension while stalk cell proliferation provides necessary material to support the extension.

\subsection{Capillary regression\label{ang_3}}
\begin{figure}[htbp]
\begin{center}
\psfig{file=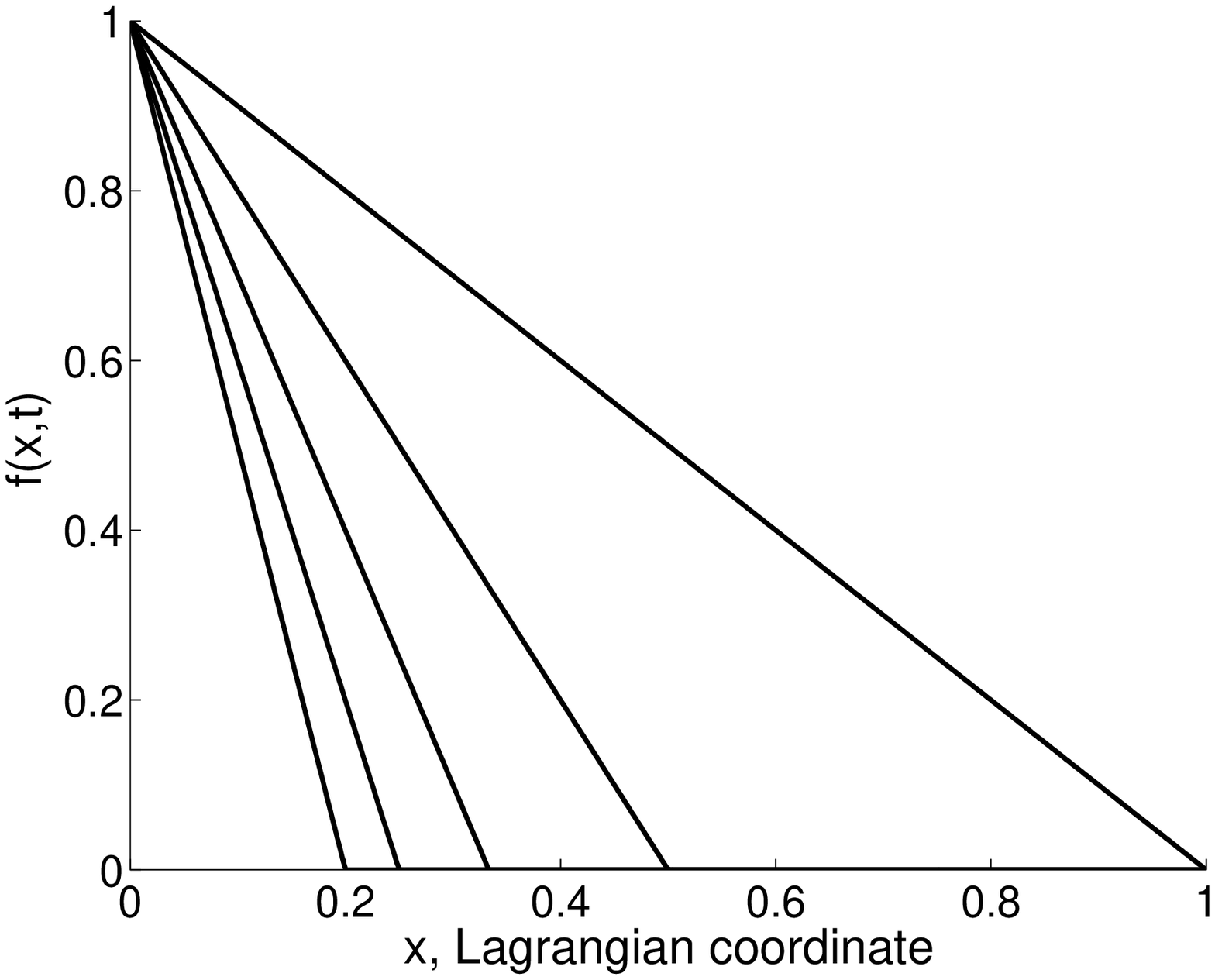, width=2.5in}[a]
\psfig{file=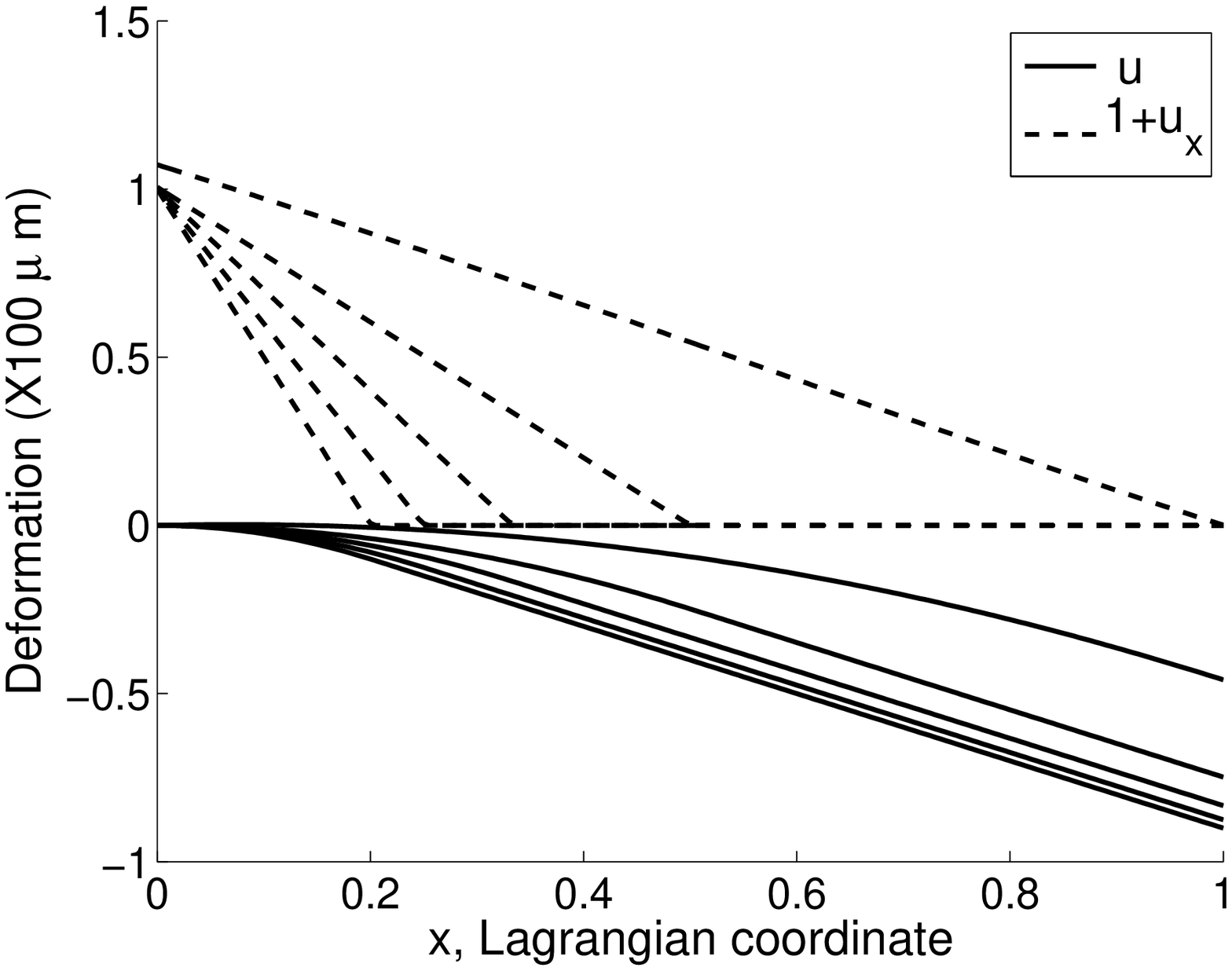, width=2.5in}[b]
\caption{Capillary regression. [a] Cell density $f$ along the capillary at time $t=1,2,3,4,5$ (from upper to lower). [b] Deformation $u$  curves (solid) and  $1+u_x$  curves (dashed) are at $t=1,2,3,4,5$, both from upper to lower.}
\label{fig3}
\end{center}
\end{figure}

The newly formed blood vessels in pathological corneal angiogenesis can be treated by the drug Bevacizumab as in \cite{Dastjerdi04012009} to induce apoptosis of endothelial cells, thus the regression of blood vessels. To simulate this therapy, we reset the friction as a constant $\beta=1$, and the cell density  as
$f=\max(1-tx, 0)$ (see Fig\,\ref{fig3}[a]), which represents the cells begin to die from the tip to the root. Because the tip cell is dying, it loses the ability to produce the protrusion force. Therefore, $g=0$.  The numerical results in Fig.\ref{fig3}[b] show that the deformation of the capillary tip is $-0.9$ at $t=5$, which means the capillary shrinks by $90\,\mu m$ after 5 days' regression, so only $10\,\mu m$ is left. The deformation gradient $1+u_x$ at $t=5$ is zero for $x\in [0.2,1]$, indicating the cells initially in this interval have died out. This simulation basically captures the experiments observed in human corneal angiogenesis treatment results in \cite{Dastjerdi04012009}.

Contrary to the last example, the cell density decreases with respect to time in capillary regression.  
In this model, the stress/stain is supported by the density. Therefore, the loss of density induces the decrease of stain, thus the regression. This phenomenon is almost impossible to be modeled by reaction-convection-diffusion models because they lack the  mechanisms to retract the capillary. For example, in the model of \cite{Anderson} where only the tip cell is modeled, the tip cell would stay where it used to be in space once the chemotactic cues are removed by therapy, thus regression observed in experiments would not occur.

\section{Conclusion \label{conclusion}}
This article analyzes a one-dimensional blood vessel growth model with both analytical and numerical tools, and demonstrates its power in simulating various biological growth cases.

Theorem\,\ref{thmexistence} concludes that the global biological solution exist if (a) the cell mass density is a decreasing function along the capillary from the root  to the tip,  or (b) the amplitude of the cell density along the capillary is less than the unity. The system may break down if these two conditions are violated. Indeed numerical examples indicate that if there is a big pressure drop from the tip to the root direction or a big oscillating amplitude of density, then the solution would probably blow up at the low density point. The underlying mechanism of breakdown is the failure of linear elasticity to resist the high pressure impact. 

Besides modeling the biological growth, this current model can also be used to describe
the slow deformation of one-dimensional thermoelastic bar, where the stress tensor is the same as \eqref{eq2} except the cell density is interpreted by temperature  ({\em cf.}\,\cite{Pabst}). Therefore the analysis can be applied in the thermoelastic case. For example, if the temperature decreases from the root to the tip of the bar, or has small oscillations along it, then the global solution exists.  
However,  if the tip of the bar is overheated or a large temperature oscillation is applied on the bar, then the bar may be broken by the expansion pressure produced by temperature differentials.

The reaction-convection-diffusion models are predominant in continuous models of angiogenesis, but they encounter great difficulties in simulating capillary growth without proliferation and capillary regression. The difficulties come from the fact that  they are unable to build connections between tip cell protrusion, change of stalk cell density,  alteration of the stress, and capillary extension. In contrast, the model of \cite{ZhengJacksonKoh} regards the whole capillary as one viscoelastic material dragged by the tip cell and directly relates cell density with elastic deformation. Plus the nonlinear interactions with surroundings, this model is capable of simulating various capillary growth patterns.

\appendix
\section{Mathematical model, non-dimensionalization, and parameters}\label{secapp}
The blood vessel capillary extension is described in the reference coordinate $x\in [0, L]$ by
\begin{eqnarray}
 \frac{2\tilde\beta |1+u_x|}{r} \dot{u}&=& \frac{\partial \tilde\sigma}{\partial x},
  \label{mol2} \\
 \tilde\sigma &=& \tilde E \frac{\partial u}{\partial x} + \tilde\mu \frac{\partial \dot{u}} {\partial x } - \tilde E \frac { (\tilde F- \rho_0)}{\rho_0},  \label{mol3}\\
 u(0,t) &= & 0, \quad \sigma (L,t) = \tilde G, \label{BC4_1}\\
 u(x,0) &=& 0, \label{BC1_2} \\
\end{eqnarray}
where $L$ is the initial capillary length, $u(x,t)$ is the displacement of material point $x$ at time $t$, 
$\tilde \beta$ is the friction, $r$ is the capillary radius,
$\tilde E$ is the Young's modulus of endothelial cells, $\tilde \mu$ is the dynamic viscosity, $\rho_0$ is the initial cell density, $\tilde F$ is the cell density multiplied with  $|1+u_x|$, and $\tilde G$ is the protrusion force at the tip of the capillary. For the derivation of this model and other  related models, see \cite{ZhengJacksonKoh}.

Choose $L$ as the characteristic length, which is $100 \,\mu m$  based on the experiments in \cite{Sholley}. Because angiogenesis typically runs a couple of days from the initiation to the penetration of tumor/tissue, we choose $T=1\,day$ as the characteristic time scale. Let $x=x'L, u=u'L, t=t' T$, and insert these relations to the above system, where we then delete the primes for simplicity. Therefore we obtain for $x \in [0, 1]$,
\begin{eqnarray}
 \beta |1+u_x| \dot{u}&=& \frac{\partial \sigma}{\partial x}, \label{non3} \\
 \sigma &=& \frac{\partial u}{\partial x} + \mu \frac{\partial \dot{u}} {\partial x} - (f - 1),
 \label{non4}\\
 u(0,t) &=& 0, \quad \sigma (1,t) = g, \label{non8}\\
 u(x,0) &=& 0, \label{non7} 
\end{eqnarray}
where $\beta =  \frac{2L^2 \tilde \beta }{r\tilde ET}$,  $\mu = \frac {\tilde\mu}{\tilde ET}$, $f=\frac{\tilde F}{\rho_0}$, and $g=\frac{\tilde G}{\tilde E}$.

 The Young's modulus $\tilde E$ for endothelial cells is between  $1.5\sim 5.6 \times 10^{3} \,\frac{pN}{\mu m^2}$ according to \cite{Costa}. The viscosity $\tilde\mu$ is not available for endothelial cells, so we replace it with the value for  fibroblasts  \cite{Ott}: $\tilde\mu=10^4\,\frac{pN \cdot s}{\mu m^2}$ (comparing with   $10^{-3}\,\frac{pN \cdot s}{\mu m^2}$ of water, $3\times 10^4\,\frac{pN \cdot s}{\mu m^2}$ of tar,  and $2.30\times 10^8 \frac{pN \cdot s}{\mu m^2}$ of pitch , all at $20^\circ C$ \cite{CRC}).  The estimate of $\tilde\beta$ will be highly dependent on the cell-environment contacts.  For example, in the case of bovine aortic endothelial cells  (BAECs) spreading on polyacrylamide gels, the friction is deduced to be $\tilde\beta=10^3\,\frac{pN \cdot s}{{\mu m}^3}$ \cite{Larripa}. However, ECs experience much stronger friction or resistance in the {\em in vivo} condition, because ECs are in association with surrounding cells such as pericytes. Therefore, we choose a large range for the friciton: $2.5\times 10^3\sim 6.5\times 10^6\,\frac{pN \cdot s}{{\mu m}^3}$. The radius of blood vessel capillary is about $10\,\mu m$. The protrusion force $\tilde G$  is about $10^{4}\,\frac{pN}{\mu m^2}$ as measured in \cite{Prass}.  With these values, the non-dimensionalized parameters become $\mu = 10^{-4}$, $\beta\in [0.01, 100]$,  and $g\in [1.7, 6.7]$.

Both the viscosity and friction resist the motion of capillary, and their relationship  can be analyzed from a model problem. Assume the deformation $u(x,t)$ is a periodic function in the free space with period $2\pi$, satisfying $u(0,t)=u(2\pi,t)=0$, and
\begin{equation}
 \beta u_t = u_{xx} + \mu u_{txx} \label{model1}.
\end{equation}
Express $u(x,t)$ in the Fourier series $\sum\limits_{k=1}^{\infty} \hat{u}_k(t) \sin(k x)$, and insert it to \eqref{model1}. Then we obtain
\begin{equation}
 (\beta+\mu k^2) \hat{u}_{k,t} = - k^2\hat{u}_k \label{model2}
\end{equation}
for each mode number $k\ge 1$.
The time scale for each mode to reach the steady state is $\beta/k^2 + \mu$. Therefore, the time scale for the whole system is $\beta + \mu$.
Given the above values of $\mu$ and $\beta$, the viscosity is far less than the friction. This implies the viscosity would have negligible effects compared with the friction on the capillary extension.

{\bf Acknowledgement.}
We thank Dapeng Du in Northeast Normal University in China for insightful discussions. Xie thanks  Jeffrey Rauch in University of Michigan for helpful discussions. Part of this work was done when Xie was visiting the Institute of Mathematical Sciences, The Chinese University of Hong Kong. He thanks the institute for its support and hospitality.

\end{document}